\documentclass[12pt]{amsart}
\usepackage{amsmath,amsthm,amssymb}
\usepackage{graphicx}
\begin{document}
\newtheorem{thm}{Theorem}
\newtheorem{pro}[thm]{Proposition}
\newtheorem{cor}[thm]{Corollary}
\newtheorem{lem}[thm]{Lemma}
\newtheorem{dfn}[thm]{Definition}
\newtheorem{rem}[thm]{Remark}
\newtheorem{prob}[thm]{Problem}
\newtheorem{exam}[thm]{Example}
\newtheorem{conj}[thm]{Conjecture}
\renewcommand{\theequation}{\arabic{section}.\arabic{equation}}
\renewcommand{\labelenumi}{\rm{(\arabic{enumi})}}
\title[Symplectic topology of Lagrangian submanifolds of ${\mathbb C}P^n$]
{Symplectic topology of Lagrangian submanifolds of ${\mathbb C}P^n$
with intermediate minimal Maslov numbers}
\author[H. Iriyeh]
{Hiroshi Iriyeh}
\date{}

\begin{abstract}
We examine symplectic topological features of certain family of monotone
Lagrangian submanifolds in ${\mathbb C}P^n$.
Firstly, we give a cohomological restriction for Lagrangian submanifolds in ${\mathbb C}P^n$
whose first integral homologies are $3$-torsion.
In particular, in the case where $n=5,8$, 
we prove the cohomologies with coefficients in $\mathbb Z_2$ of such Lagrangian submanifolds
are isomorphic to that of $SU(3)/(SO(3) {\mathbb Z}_3)$ and $SU(3)/{\mathbb Z}_3$, respectively.
Secondly, we calculate the Floer cohomology of a monotone Lagrangian submanifold $SU(p)/{\mathbb Z}_p$
in ${\mathbb C}P^{p^2-1}$ with coefficients in $\mathbb Z_2$ by using Biran-Cornea's theory.
\end{abstract}

\maketitle

\section{Introduction and main results}

Let $(M,\omega)$ be a symplectic manifold,
i.e., $M$ is a smooth manifold with a closed nondegenerate two-form $\omega$.
A submanifold $L$ of $M$ is called Lagrangian if
$\dim_{\mathbb R}L=(1/2)\dim_{\mathbb R}M$ and the restriction of $\omega$ on $L$ vanishes.
Throughout this paper all symplectic manifolds are assumed to be {\it tame}, i.e.,
there exists an almost complex structure $J$ such that the bilinear form $\omega(\cdot,J\cdot)$ defines
a Riemannian metric on M which is geometrically bounded (see \cite{ALP94}).
And all Lagrangian submanifolds are assumed to be closed (i.e., compact and without boundary),
connected and embedded.

In the complex projective space ${\mathbb C}P^n$,
there are two familiar examples of Lagrangian submanifolds.
One is the real form ${\mathbb R}P^n$ of it and the other is the Clifford torus defined by
${\mathbb T}_{\rm clif}^n=\{ [z_0:\cdots:z_n] \in {\mathbb C}P^n \mid |z_0|=\cdots=|z_n| \}$.
For instance, the Arnold-Givental conjecture was first proved for the former example (see \cite{Oh93'})
from the view point of Floer theory,
and then researches of general cases were developed
(for instance, \cite{Oh95}, \cite{Frauenfelder04}, \cite{FOOO} and \cite{IST13}).
On the other hand, the latter example is a maximal orbit of the standard $n$-dimensional torus action on ${\mathbb C}P^n$.
The calculation of the Floer cohomology of ${\mathbb T}_{\rm clif}^n$ was carried out by C.\ H.\ Cho \cite{Cho04} and,
at present, that of general Lagrangian torus orbits in toric Fano manifolds has been intensively studied.
However, it seems that, until now, there is no Floer theoretic 
study of Lagrangian submanifolds in ${\mathbb C}P^n$ beyond ${\mathbb R}P^n$ and ${\mathbb T}_{\rm clif}^n$.

The purpose of the present paper is to initiate systematic symplectic topological research of a certain class of
monotone Lagrangian submanifolds in ${\mathbb C}P^n$ naturally including the above two examples.
In particular, we present new ideas to obtain results about {\it homological rigidity {\rm (or} uniqueness})
and {\it non-displaceability} of such Lagrangian submanifolds.

We now recall the definitions of monotoneness and the minimal Maslov number of a Lagrangian submanifold.
For a Lagrangian submanifold $L$ in a symplectic manifold $(M,\omega)$, two homomorphisms
\begin{equation*}
I_{\mu,L} : \pi_2(M,L) \to {\mathbb Z},\ \ 
I_\omega : \pi_2(M,L) \to {\mathbb R}
\end{equation*}
are defined as follows.
For a smooth map $w : (D^2,\partial D^2) \to (M,L)$,
$I_{\mu,L}(w)$ is defined to be the Maslov number of the bundle pair
$(w^*TM,$ $(w|\partial D^2)^*TL)$ 
and $I_\omega$ is defined by $I_\omega(w)=\int_{D^2}w^*\omega$.
These two homomorphisms do not depend on the representative $w$.
Then $L$ is said to be {\it monotone} if there exists a constant $\alpha>0$
such that $I_\omega = \alpha I_{\mu,L}$.
The {\it minimal Maslov number} $N_L$ of $L$ is defined to be
the positive generator of $\mathrm{Im}(I_{\mu,L}) \subset {\mathbb Z}$.

The present paper mainly concentrates on the case of the complex projective space $({\mathbb C}P^n,\omega_{\rm FS})$
with the standard Fubini-Study K\"ahler form $\omega_{\rm FS}$.
First of all, we shall provide certain family of homogeneous Lagrangian submanifolds in ${\mathbb C}P^n$
which are monotone. 

\begin{pro} \label{pro:A-O}
There exist following four families of monotone Lagrangian submanifolds $L$ in ${\mathbb C}P^n:$
\begin{enumerate}
\item $\displaystyle L=\frac{SU(p)}{\mathbb Z_p}$ \quad if \quad $n+1=p^2$,
\item $\displaystyle L=\frac{SU(p)}{SO(p){\mathbb Z}_p}$ \quad if \quad $\displaystyle n+1=\frac{p(p+1)}{2}$,
\item $\displaystyle L=\frac{SU(2p)}{Sp(p){\mathbb Z}_{2p}}$ \quad if \quad $n+1=p(2p-1)$,
\item $\displaystyle L=\frac{E_6}{F_4{\mathbb Z}_3}$ \quad if \quad $n+1=27$,
\end{enumerate}
where $p \in \mathbb N \setminus \{ 1 \}$.
The minimal Maslov numbers of $L$ in the first and second families equal $2p$ and $p+1$, respectively,
if $p$ is prime.
That of the example $(4)$ equals $18$.
\end{pro}

Here we observe the following inequality concerning the minimal Maslov number $N_L$.

\begin{pro} \label{pro:1--n+1}
Let $L$ be a monotone Lagrangian submanifold of ${\mathbb C}P^n$.
Then we have
\begin{equation} \label{eq:Maslovbound}
1 \leq N_L \leq n+1,
\end{equation}
and if $N_L=n+1$ holds, then $L$ is a $\mathbb Z_2$-homological ${\mathbb R}P^n$.
\end{pro}

This is an easy consequence of the Lagrangian circle bundle construction \cite[Proposition 4.1.A]{Biran06}
and a result by Oh (see Theorem \ref{thm:Oh96} below) which states that $1 \leq N_L \leq {n+1}$
for any monotone Lagrangian submanifold $L \subset {\mathbb C}^{n+1}$.
The upper equality condition of (\ref{eq:Maslovbound}) rephrases a result by P. Biran
(see Theorem \ref{thm:Biran1} and Remark \ref{rem:Biran1} below).
We note that $N_{{\mathbb T}_{\rm clif}^n}=2$, which is the minimal value of $N_L$ among orientable 
Lagrangian submanifolds of ${\mathbb C}P^n$.
Thus each Lagrangian submanifold in the list in Proposition \ref{pro:A-O} possesses
an {\it intermediate} minimal Maslov number.

In the following subsections we explain some symplectic topological properties, e.g.,
homological rigidity, non-displaceability and uniruling, of the Lagrangian submanifolds in Proposition \ref{pro:A-O}.
We believe that their model Lagrangian submanifolds also provide interesting examples for other aspects of symplectic topology,
for instance, Lagrangian quantum homology, symplectic quasi-state, etc.\
as well as contribute to the classification of monotone Lagrangian submanifolds in ${\mathbb C}P^n$.

\subsection{Homological rigidity of Lagrangian submanifolds in ${\mathbb C}P^n$}

In this subsection we shall treat {\it homological rigidity} or {\it uniqueness} of Lagrangian submanifolds.
This phenomenon means low-dimensional topological invariants of Lagrangian submanifolds
determine their entire (co)homology.
For a closed symplectic manifold $M$,
this phenomenon was first discovered by P.\,Seidel \cite{Seidel00} for the case $M={\mathbb C}P^n$.
Nowadays, many results concerning homological rigidity are known including other symplectic manifolds
(see \cite{Biran06}, \cite{Buhovsky04}, \cite{Damian12} and references therein).
Here we review results by P. Biran and M. Damian which are the starting point of the present research
(see also \cite[Corollary 1.2.11]{Biran-Cornea09}).

\begin{thm}[Biran \cite{Biran06}, Theorem A] \label{thm:Biran1}
Let $L$ be a Lagrangian submanifold of ${\mathbb C}P^n$ such that $H_1(L;\mathbb Z)$ is $2$-torsion,
i.e., $2 H_1(L;\mathbb Z)=0$.
Then
\begin{enumerate}
\item $H^*(L;{\mathbb Z}_2) \cong H^*({\mathbb R}P^n;{\mathbb Z}_2)$ as graded vector spaces,

\item if $n$ is even, then the isomorphism in $(1)$ is an isomorphism as graded algebras.
\end{enumerate}
\end{thm}

\begin{thm}[Damian \cite{Damian12}, Theorem 1.8.c] \label{thm:Damian1}
Under the same assumption as Theorem \ref{thm:Biran1},
if $n$ is odd, then the Lagrangian submanifold $L \subset {\mathbb C}P^n$ satisfies $\pi_1(L) \cong \mathbb Z_2$
and the universal cover of $L$ is homeomorphic to $S^n$.
\end{thm}

\begin{rem} \label{rem:Biran1} \rm
The assumption $2 H_1(L;\mathbb Z)=0$ implies that $L \subset {\mathbb C}P^n$ is monotone and $N_L=n+1$.
Biran's proof \cite[p.\ 313]{Biran06} of Theorem \ref{thm:Biran1} actually shows that
these two conditions yield the conclusion.
Hence, Theorem \ref{thm:Biran1} can be considered as a homological characterization of monotone Lagrangian submanifolds
$L \subset {\mathbb C}P^n$ with the maximal $N_L$ (see Proposition \ref{pro:1--n+1}).
\end{rem}

Now we present a new class of homologically rigid Lagrangian submanifolds beyond the case of ${\mathbb R}P^n$.
The following is the main result.

\begin{thm} \label{thm:main1}
Let $L$ be a Lagrangian submanifold of $\mathbb CP^n$ such that $H_1(L;\mathbb Z)$ is $3$-torsion,
i.e., $3 H_1(L;\mathbb Z)=0$.
Then $L$ is monotone and orientable, and $3\,|\,n+1$ and $n \geq 5$ hold.
Moreover
\begin{enumerate}
\item If $n=5$, then
$$
H^*(L; \mathbb Z_2) \cong H^*\left(\frac{SU(3)}{SO(3) \mathbb Z_3}; \mathbb Z_2\right)
$$
as graded vector spaces.
\item If $n=8$, then
$$
H^*(L; \mathbb Z_2) \cong H^*\left(\frac{SU(3)}{\mathbb Z_3}; \mathbb Z_2\right)
$$
as graded algebras.
\item If $n=26$ and $H^i(L; \mathbb Z_2) \cong 0 \ (i=2,3,4)$, then
$$
H^*(L; \mathbb Z_2) \cong H^*\left(\frac{E_6}{F_4 \mathbb Z_3}; \mathbb Z_2\right)
$$
as graded algebras.
\end{enumerate}
Furthermore, the Euler characteristic $\chi(L)=\sum_i (-1)^i \dim_{\mathbb Z_2}H_i(L;\mathbb Z_2)$
of $L$ is equal to zero.
\end{thm}

\begin{rem} \label{rem:coh} \rm
In the following, we denote by $\wedge_{\mathbb Z}$ and $\wedge_2$ exterior algebra over $\mathbb Z$ and $\mathbb Z_2$,
respectively, with generators in round brackets.
It is known that
\begin{eqnarray*}
& & H^*\left(\frac{SU(3)}{SO(3) {\mathbb Z}_3};{\mathbb Z}_2\right) \cong \wedge_2(x_2,x_3), \
 H^*\left(\frac{SU(3)}{{\mathbb Z}_3};{\mathbb Z}_2\right) \cong \wedge_2(x_3,x_5), \\
& & H^*\left(\frac{E_6}{F_4{\mathbb Z}_3};{\mathbb Z}_2\right) \cong \wedge_2(x_9,x_{17})
\end{eqnarray*}
as graded algebras (see Section 2 for details).
Here the suffix $i$ of $x_i$ denotes the degree of $x_i$.
\end{rem}

\subsection{Non-displaceability of model Lagrangian submanifolds}

\
Next we turn to a problem whether model Lagrangians of ${\mathbb C}P^n$ are displaceable or not.
A diffeomorphism $\phi$ of $(M,\omega)$ is called {\it Hamiltonian} if it is the time-one map of the flow which is
defined by a compactly supported time dependent Hamiltonian function on $M$.
A Lagrangian submanifold $L \subset M$ is said to be {\it displaceable} if
there exists a Hamiltonian diffeomorphism $\phi \in \mathrm{Ham}(M,\omega)$ such that
$$
L \cap \phi(L) = \emptyset,
$$
otherwise, $L$ {\it non-displaceable}.
If the Floer cohomology $HF(L,L)$ with coefficients in $\mathbb Z_2$ does not vanish,
then we see that $L$ is non-displaceable in $M$.
Specifically, our method is well applicable to the case
$L=SU(p)/\mathbb Z_p \subset {\mathbb C}P^{p^2-1}$.

\begin{thm} \label{thm:non-displaceable4}
Let $L$ be the Lagrangian submanifold $SU(p)/{\mathbb Z}_p \subset {\mathbb C}P^{p^2-1}$,
where $p$ is a power of $2$.
Then the Floer cohomology $HF(L,L)$ with coefficients in $\mathbb Z_2$ is isomorphic to
$H^*(L;{\mathbb Z}_2) \otimes \Lambda$, where $\Lambda={\mathbb Z}_2[T,T^{-1}]$.
In particular, $L \subset {\mathbb C}P^{p^2-1}$ is non-displaceable.
\end{thm}

In contrast, we have  $HF(L,L)=0$ for an odd number $p$ (see Corollary \ref{cor:vanishingHF}).
Nevertheless, using the Floer cohomology with coefficients in $\mathbb Z$, we prove

\begin{pro} \label{pro:non-displaceable}
The Lagrangian submanifold $SU(3)/{\mathbb Z}_3 \subset {\mathbb C}P^8$ is non-displaceable.
\end{pro}

Notice that $\chi(SU(p)/{\mathbb Z}_p)=0$ for any $p \in {\mathbb N} \setminus \{ 1 \}$.


\subsection{Uniruling of model Lagrangian submanifolds in ${\mathbb C}P^n$}

For some model Lagrangian submanifolds $L$, we can prove the existence of a pseudo-holomorphic disc
with its boundary on $L$.
We review here a combinient terminology for such results (see \cite[Definition 1.1.2]{Biran-Cornea09}).

\begin{dfn} \rm
A Lagrangian submanifold $L \subset (M,\omega)$ is said to be {\it uniruled of order $k$}
if for any point $Q \in L$, there exists a generic family of almost complex structures $\mathcal{J}$
with the property that for each $J \in \mathcal{J}$ there exists a nonconstant $J$-holomorphic disc
$u:(D^2,\partial D^2) \to (M,L)$ such that
$$
Q \in u(\partial D^2) \quad {\rm and} \quad \mu(u) \leq k,
$$
where $\mu(u)=I_{\mu,L}(u)$ is the Maslov number of $u$.
\end{dfn}

Combining a result by Biran and Cornea (see Theorem \ref{thm:Biran-Cornea1}) with Corollary \ref{cor:vanishingHF}
in Section 5, we prove

\begin{cor} \label{cor:uniruling}
The monotone Lagrangian submanifold $SU(p)/{\mathbb Z}_p \subset {\mathbb C}P^{p^2-1}$, where $p$ is an odd prime,
is uniruled of order $2p$.
\end{cor}

\subsection{Organization of the paper}

The rest of the paper is organized as follows.
In Section 2, we give some preliminary properties of model Lagrangian submanifolds,
in particular, their minimal Maslov numbers.
In Section 3, we briefly recall the definitions of Floer cohomology for Lagrangian submanifolds
and the spectral sequence introduced by Oh and developed by Biran,
including its multiplicative structure established by L.\ Buhovsky.
Section 4 is devorted to proving the main result (Theorem \ref{thm:main1}).
To prove (2) and (3) in the theorem we use, instead of the Lagrangian circle bundle construction,
a recent result by Biran and Khanevsky (see Theorem \ref{thm:B-K}).
It relates the Floer cohomology $HF(L,L)$ with coefficients in $\mathbb Z_2$ with
the $\mathbb Z_2$-Euler class of the normal bundle $\mathcal{N}|_L$ defined by a subcritical polarization.
In the present case, the polarization is $({\mathbb C}P^{n+1},\omega_{\rm FS},J,{\mathbb C}P^n)$ and
we can prove the vanishing of the Euler class of $\mathcal{N}|_L$ under a suitable condition,
and hence we obtain $HF(L,L)=0$ under a further additional condition, 
which is a new approach we introduce in the present paper.
By using the spectral sequence, we obtain topological constraints on Lagrangian submanifolds in ${\mathbb C}P^n$.
In Section 5, we prove some new results about the non-displaceability of certain model Lagrangian submanifolds
in ${\mathbb C}P^n$ by two approaches.
One is a use of Biran and Cornea's result (see Theorem \ref{thm:Biran-Cornea1}).
The other is a use of the spectral sequence and the Floer cohomology with coefficients in $\mathbb Z$.
In Section 6, as an application, we show the existence of a pseudo-holomorphic disc with its boundary on
a model Lagrangian submanifold.
Such Lagrangian turns out to be a border case of Biran and Cornea's theorem mentioned above.
This fact together with results in the previous section yields certain uniruling results for
such model Lagrangian submanifolds.
Finally, in Section 7 we notice that results in this paper will provide us with a step toward
the classification of monotone Lagrangian submanifolds in ${\mathbb C}P^n$.

\section{Model Lagrangian submanifolds}

Let us consider again the Lagrangian submanifolds in ${\mathbb C}P^n$ introduced in the previous section:
$$
\frac{SU(p)}{\mathbb Z_p}, \quad \frac{SU(p)}{SO(p){\mathbb Z}_p}, \quad
\frac{SU(2p)}{Sp(p){\mathbb Z}_{2p}}, \quad \frac{E_6}{F_4{\mathbb Z}_3} \subset {\mathbb C}P^n,
$$
where $n$ is an appropriate integer in Proposition \ref{pro:A-O}.
These examples first appeared in the paper \cite{Amarzaya-Ohnita03} by Amarzaya and Ohnita
in the context of research on minimal submanifold theory.
Note that R. Chiang \cite{Chiang04} rediscovered the first family in the above list
from the view point of momentum map.
Including ${\mathbb R}P^n$, they are all irreducible {\it embedded} minimal submanifolds in ${\mathbb C}P^n$
(i.e., critical points of volume functional of $L$ with respect to the induced Riemannian metric
from ${\mathbb C}P^n$).
Indeed, they possess parallel second fundamental forms,
and hence they are not only homogeneous spaces but also symmetric spaces.
Although the Clifford torus ${\mathbb T}_{\rm clif}^n$ is not irreducible, it also has such a property.
The local classification of symmetric Lagrangian submanifolds of ${\mathbb C}P^n$ was given by
Naitoh and Takeuchi (see \cite[Theorem 4.5]{Naitoh81}, \cite{Naitoh-Takeuchi82}).
Moreover the cohomology rings of $SU(p), SU(p)/SO(p)$ and $E_6/F_4$,
i.e., universal covering spaces of the model spaces, are known
(see, e.g., \cite{Borel54} and \cite[Proposition 2.5]{Araki61})\,:
\begin{eqnarray*}
& & H^*(SU(p);\mathbb Z) \cong \wedge_{\mathbb Z}(x_3,x_5,\ldots,x_{2p-1}), \\
& & H^*\left(\frac{SU(p)}{SO(p)};{\mathbb Z}_2\right) \cong \wedge_2(x_2,x_3,\ldots,x_p), \\
& & H^*\left(\frac{E_6}{F_4};{\mathbb Z}\right) \cong \wedge_{\mathbb Z}(x_9,x_{17}).
\end{eqnarray*}
Let us also treat the case of model Lagrangian submanifolds with coefficients in $\mathbb Z_2$.
For instance, consider the case $SU(p)/{\mathbb Z}_p$, where $p$ is an odd number greater than or equal to three.
Since it has no $2$-torsion in its homology, we have
\begin{eqnarray} \label{eq:cohSU}
\ \ \ \ \ \ \ \ H^*\left(\frac{SU(p)}{\mathbb Z_p};\mathbb Z_2\right) \cong H^*(SU(p);\mathbb Z_2)
\cong \wedge_2(x_3,x_5,\ldots,x_{2p-1}).
\end{eqnarray}
We also have
$$
H^*\left(\frac{SU(3)}{SO(3) {\mathbb Z}_3};{\mathbb Z}_2\right) \cong \wedge_2(x_2,x_3),\ 
H^*\left(\frac{E_6}{F_4{\mathbb Z}_3};{\mathbb Z}_2\right) \cong \wedge_2(x_9,x_{17}).
$$

\smallskip

It is well-known that ${\mathbb R}P^n \subset {\mathbb C}P^n$ is monotone and
$N_{{\mathbb R}P^n}=n+1$.
Similarly we can prove Proposition \ref{pro:A-O}.

\smallskip

\noindent
{\bf Proof of Proposition \ref{pro:A-O}.} \
It is easy to check that a Lagrangian submanifold $L$ of ${\mathbb C}P^n$ satisfying $\pi_1(L) \cong \mathbb Z_p$
is monotone (see Lemma \ref{lem:minimal Maslov} below for a similar argument).

Assume that $p$ is prime.
We define two subgroups $\Gamma_\omega, \Gamma_{\omega,L}$ of $\mathbb R$ by
$$
\Gamma_\omega=\{[\omega](A) \mid A \in H_2(M; \mathbb Z) \}, \ \
\Gamma_{\omega,L}=\{[\omega](B) \mid B \in H_2(M,L; \mathbb Z) \}.
$$
Recall that a symplectic manifold $(M,\omega)$ is said to be {\it prequantizable} if $\Gamma_\omega$ is either
trivial or discrete.
A Lagrangian submanifold $L \subset M$ is called {\it cyclic} if $\Gamma_{\omega,L} \subset \mathbb R$ is discrete.
Note that $\Gamma_\omega$ is a subgroup of $\Gamma_{\omega,L}$.
When $L$ is cyclic, we can define a positive integer $n_L:=\#(\Gamma_{\omega,L}/\Gamma_\omega)$
(see \cite[p.\ 473]{Oh94}).
It is well-known (see, for instance, \cite{Kostant70}) that if $(M,\omega)$ is prequantizable, then there is a principal
$\mathbb R/\Gamma_\omega$-bundle $\pi:Q \to M$ with a connection form $\theta$ such that its curvature form $d\theta$
satisfies
$$
d\theta = \pi^* \omega.
$$
Hence, given any Lagrangian submanifold $L \subset ({\mathbb C}P^n,\omega_{\rm FS})$ the connection $\theta|_L$
on $Q|_L$ is flat.
Then
$$
G_L(x):=\mathrm{im} \{ \pi_1(L,x) \to \mathbb R/\Gamma_\omega \cong S^1 \}
$$
is nothing but the holonomy group of $Q|_L$ with a base point $x \in L$.
Since $\pi_1(L) \cong \mathbb Z_p$ and $p$ is prime, the group $G_L(x)$ is isomorphic to either $\mathbb Z_p$ or trivial.
Here we use the following formula (see \cite[Proposition 3.6]{Ono05}):
\begin{equation} \label{eq:Ono}
2(n+1)=n_L N_L.
\end{equation}
If the latter case happens, then $n_L=1$ and $N_L=2(n+1)$ hold.
But, it is impossible from (\ref{eq:Maslovbound}).
Hence we have $G_L(x) \cong \mathbb Z_p$, that is, $n_L=p$.
Thus we obtain $N_L=2(n+1)/p$.
This immediately yields the latter half of the statement.
\hfill \qed

\begin{rem} \rm 
Let $L$ be one of the Lagrangian submanifolds in the family (1) or (2) in Proposition \ref{pro:A-O}.
When $p$ is {\it not} prime, we have at least an inequality
\begin{equation} \label{eq:MMineq}
N_L \geq \frac{2(n+1)}{p}.
\end{equation}
Indeed, since the order $n_L$ of the group
$G_L(x)=\mathrm{im} \{ \pi_1(L,x) \cong \mathbb Z_p \to S^1 \}$
is at most $p$, we have (\ref{eq:MMineq}) from (\ref{eq:Ono}).
\end{rem}

\section{Floer cohomology and the quantum cup product on it}

In this section, we briefly review the Lagrangian Floer theory
as developed in \cite{Oh93}, \cite{Oh96}, \cite{Biran06} and \cite{Buhovsky10}.

\subsection{Basic settings}

Let $(M,\omega)$ be a tame symplectic manifold.
Let $L_0$ and $L_1$ be closed monotone Lagrangian submanifolds of $M$ which intersect transversally.
Assume that $N_{L_i} \geq 3$ for $i=0,1$.
We choose a time-dependent family $J= \{ J_t \}_{0 \leq t \leq 1}$  
of almost complex structures on $M$ compatible with
the symplectic form $\omega$.
The Floer chain complex $CF(L_0,L_1)$ is the vector space over
$\mathbb Z_2$ generated by the finitely many elements of $L_0 \cap L_1$.
For a generic choice of $J$, the boundary operator $d_J : CF(L_0,L_1) \to CF(L_0,L_1)$ is defined by
counting $J$-holomorphic strips in $M$ connecting pairs of points of $L_0 \cap L_1$ and
with boundary in $L_0$ and $L_1$.
In this setting, we have

\smallskip

(i) \ $d_J \circ d_J=0$, so the cohomology $H^*(CF(L_0,L_1),d_J)$ is well-defined;
it denoted by $HF(L_0,L_1; J)$ and called the Floer cohomology of the pair $(L_0,L_1)$;

(ii) \ $HF(L_0,L_1; J)$ is independent of the choice of $J$;

(iii) \ if $L_1$ and $L_1'$ are Hamiltonian isotopic, then $HF(L_0,L_1) \cong HF(L_0,L_1')$.

\smallskip

These foundational results were established by Floer \cite{Floer88} and Oh \cite[Theorems 4.4, 4.5]{Oh93}.
In the present paper, we use $\mathbb Z_2$ as coefficients of the Floer cohomology, except for
the proof of Proposition \ref{thm:main2}.

From now on, we consider the case where $L_1$ is Hamiltonian isotopic to $L:=L_0$, i.e.,
$L_1=\phi(L),\ \phi \in \mathrm{Ham}(M,\omega)$.
In this setting, the Floer cohomology is also well-defined for the case $N_L=2$,
and the Floer complex $CF(L,\phi(L))$ has a relative $\mathbb Z_{N_L}$-grading.
It depends on a choice of a base intersection point $x_0 \in L \cap \phi(L)$.
Another choice of such a point yields a shift in the grading.
For fixed $x_0$, we denote by $CF^{i ({\rm mod}\ N_L)}(L,\phi(L); x_0)$ the $i$-th (mod $N_L$) component of $CF(L,\phi(L))$.
It is shown that the differential $d_J$ increases grading by one, i.e.,
$d_J: CF^{* ({\rm mod}\ N_L)}(L,\phi(L); x_0) \to CF^{*+1 ({\rm mod}\ N_L)}(L,\phi(L); x_0)$
(see \cite{Oh96}).
Hence, the Floer cohomology
$$
HF(L,\phi(L); J)=\bigoplus_{i=0}^{N_L-1} HF^{i ({\rm mod}\ N_L)}(L,\phi(L); J,x_0)
$$
has a relative $\mathbb Z_{N_L}$-grading.

\subsection{The case of HF(L,L)}

Let $U(L)$ be a Weinstein neighbourhood of $L$ in $M$.
Let $L_\epsilon$ be a Hamiltonian perturbation of $L$ built in $U(L)$ using a $C^2$-small Morse function
$f: L \to \mathbb R$.
Assume that $f$ has exactly one relative minimum $x_0$.
Denote by $C_f^*$ the Morse complex of $f$ and we use $x_0$ as a base intersection point
for the Floer complex.
From now on, we write $CF(L,L):=CF(L,L_\epsilon)$ and $d:=d_J$.
It is shown that
$$
CF^{i ({\rm mod}\ N_L)}(L,L) = \bigoplus_{j \equiv i ({\rm mod}\ N_L)} C_f^j
$$
(see \cite{Oh96}).
Since $d: CF^{* ({\rm mod}\ N_L)}(L,L) \to CF^{*+1 ({\rm mod}\ N_L)}(L,L)$,
we can write $d=\sum_{j \in \mathbb Z} \partial_j$,
where $\partial_j: C_f^* \to C_f^{*+1-jN_L}$.
An index computation shows that $\partial_j=0$ for any $j<0$, and
$\partial_j=0$ for any $j> \nu$, where $\nu:=[(\dim L+1)/N_L]$.
Thus $d$ is decomposed into
$$
d=\partial_0+\partial_1+\cdots+\partial_\nu.
$$
In \cite[Section 4]{Oh96} Oh proved that, for a suitable choice of $J$ and a Riemannian metric on $L$,
the operator  $\partial_0: C_f^* \to C_f^{*+1}$ can be identified with the Morse differential.
Hence $H^*(C_f,\partial_0) \cong H^*(L;\mathbb Z_2)$.

\subsection{Oh-Biran's spectral sequence}

We shall now briefly review a spectral sequence which enables us to calculate the Floer cohomology
$HF(L,L)$ using the operators $\partial_1,\ldots,\partial_\nu$.
For the details, see \cite[Section 5]{Biran06}.

Let $\Lambda={\mathbb Z}_2[T,T^{-1}]$ be the algebra of Laurent polynomials over ${\mathbb Z}_2$
and we define the degree of $T$ to be $N_L$.
Then $\Lambda$ is decomposed as $\Lambda= \bigoplus_{i \in \mathbb Z}\Lambda^i$,
where $\Lambda^i$ is the subspace of homogeneous elements of degree $i$.
We refine the Morse complex $C_f$ as $\tilde{C}=C_f \otimes \Lambda$, i.e.,
$$
\tilde{C}^l(L,L) = \bigoplus_{k \in \mathbb Z} C_f^{l-k N_L}(L,L) \otimes \Lambda^{k N_L}
$$
for every $l \in \mathbb Z$ and the operator $d$ as $\tilde{d}:\tilde{C}^* \to \tilde{C}^{*+1}$,
$$
\tilde{d}=\partial_0 \otimes 1 + \partial_1 \otimes T +\cdots+ \partial_\nu \otimes T^\nu,
$$
where $T^i:\Lambda^* \to \Lambda^{*+i N_L}$ is a multiplication by $T^i$.
Then we see that $\tilde{d} \circ \tilde{d}=0$ and
$$
H^l(\tilde{C},\tilde{d}) \cong HF^{l ({\rm mod}\ N_L)}(L,L)
$$
for every $l \in \mathbb Z$.
Let us consider the following decreasing filtration on $\tilde{C}$:
$$
F^p \tilde{C} = \left\{ \sum x_i \otimes T^{n_i} \mid x_i \in C_f, n_i \geq p \right\}.
$$
Since $C_f^j=0$ for $j<0$, $\dim L <j$, the filtration $F^p \tilde{C}$ is bounded.
Let $\{ E_r^{p,q},d_r \}$ be the spectral sequence defined by this filtration.
Notice that this filtration is different from the one used by Oh \cite[p.\ 338]{Oh96}.

\begin{thm}[Biran \cite{Biran06}, Theorem 5.2.A] \label{thm:Oh-Biran}
The spectral sequence \\ $\{ E_r^{p,q},d_r \}$ has the following properties$:$
\begin{enumerate}
\item $E_0^{p,q} = C_f^{p+q-p N_L} \otimes \Lambda^{p N_L},\ d_0 = [\partial_0] \otimes 1$.

\item $E_1^{p,q} = H^{p+q-p N_L}(L;{\mathbb Z}_2) \otimes \Lambda^{p N_L},\ d_1 = [\partial_1] \otimes T$, where
$$
[\partial_1] : H^{p+q-p N_L}(L;{\mathbb Z}_2) \to H^{p+1+q-(p+1)N_L}(L;{\mathbb Z}_2)
$$
is induced from $\partial_1$.

\item For every $r \geq 1$, $E_r^{p,q}$ has the form $E_r^{p,q} = V_r^{p,q} \otimes \Lambda^{p N_L}$
with $d_r=\delta_r \otimes T^r$,
where each $V_r^{p,q}$ is a vector space over $\mathbb Z_2$ and $\delta_r$ is a homomorphism
$\delta_r : V_r^{p,q} \to  V_r^{p+r,q-r+1}$
defined for every $p,q$ and satisfies $\delta_r \circ \delta_r=0$.
Moreover,
$$
V_{r+1}^{p,q} = \frac{\ker(\delta_r : V_r^{p,q} \to  V_r^{p+r,q-r+1})}{\mathrm{im}(\delta_r : V_r^{p-r,q+r-1} \to  V_r^{p,q})}.
$$

\item $E_r^{p,q}$ collapses at $(\nu+1)$-step, namely $d_r=0$ for every $r \geq \nu+1$,
and the sequence converges to $HF(L,L)$, i.e.,
$$
\bigoplus_{p+q=l} E_\infty^{p,q} \cong HF^{l ({\rm mod}\ N_L)}(L,L)
$$
for every $l \in \mathbb Z$, where $\nu=[(\dim L+1)/N_L]$.

\item For all $p \in \mathbb Z$, we have
$\oplus_{q \in \mathbb Z} E_\infty^{p,q} \cong HF(L,L)$.
\end{enumerate}
\end{thm}

In this paper, we refer to $\{ E_r^{p,q},d_r \}$ as {\it Oh-Biran's spectral sequence}.

\subsection{A multiplicative structure}

The spectral sequence $\{ E_r^{p,q},d_r \}$ possesses a multiplicative structure.
In particular, the multiplication on $E_1^{p,q}$ induces that on $H^*(L;\mathbb Z_2)$
which coincides with the classical cup product on $H^*(L;\mathbb Z_2)$.
We briefly review it (see \cite{Buhovsky10} for more details).
Note that, in this paper, the product structure is used only in the proof of (2) and (3) in Theorem \ref{thm:main1}.
 
Let us consider generic Morse functions $f,g,h : L \to \mathbb R$ and denote by
$(CF_f,d^f)$, $(CF_g,d^g)$ and $(CF_h,d^h)$ the corresponding Floer complexes.
In \cite{Buhovsky10} it was defined a quantum product
$\star: CF_f \otimes CF_g \to CF_h$ which satisfies the following Leibniz rule:
$$
d^h(a \star b) = d^f(a) \star b + a \star d^g(b)
$$
for every $a \in CF_f,\,b \in CF_g$.
This product is compatible with the filtrations on $CF_f,\,CF_g$ and $CF_h$ such that
$\star: F^pCF_f \otimes F^{p'}CF_g \to F^{p+p'}CF_h$.
It induces the product on $\{ E_r^{p,q},d_r \}$ satisfying that
$\star: E_r^{p,q}(f) \otimes E_r^{p',q'}(g) \to E_r^{p+p',q+q'}(h)$
at each stage.
Then the differential $d_r$ satisfies the Leibniz rule with respect to this product,
and the product at the $(r+1)$ stage comes from the one at $r$ stage.
Moreover, these products induce products
$V_r^{p,q}(f) \otimes V_r^{p',q'}(g) \to V_r^{p+p',q+q'}(h)$
and the differential $\delta_r: V_r^{p,q} \to V_r^{p+r,q-r+1}$ satisfies the Leibniz rule.
Then the following theorem is proved.

\begin{thm}[Buhovsky \cite{Buhovsky10}, Theorem 5] \label{thm:Buhovsky10}
The product on $V_1$, induced from $\star$, coincides with the classical cup product on $H^*(L;{\mathbb Z}_2)$.
\end{thm}

\section{Homological rigidity of Lagrangian submanifolds whose $H_1$ are $3$-torsion}

In this section we prove the homological rigidity of Lagrangian submanifolds in ${\mathbb C}P^n$
whose first integral homologies are $3$-torsion stated in the introduction.
Before we discuss the case $M={\mathbb C}P^n$, let us review two well-known topological constraints
for Lagrangian submanifolds in $\mathbb C^n \cong \mathbb R^{2n}$.
They are used in the following arguments.

A Lagrangian submanifold $L \subset (M,\omega)$ is said to be {\it weakly exact} if the homomorphism
$I_\omega$ vanishes, and {\it exact} if $\omega=d\lambda$ for an one-form $\lambda$ on $M$ and
the restriction $\lambda|_L$ on $L$ is an exact one-form on $L$.
Then M.\ Gromov \cite{Gromov85} proved

\begin{thm} \label{thm:Gromov}
There is no (weakly) exact Lagrangian submanifold $L \subset \mathbb C^n$.
In particular, we have $H^1(L;\mathbb R) \neq 0$ for any Lagrangian submanifold $L \subset \mathbb C^n$.
\end{thm}

The following result about $N_L$ obtained by Oh improves previous results of L. Polterovich
\cite{Polterovich91}, \cite{Polterovich91'}.

\begin{thm}[\cite{Oh96}, Theorem 5.3] \label{thm:Oh96}
For any monotone Lagrangian submanifold $L \subset \mathbb C^n$, we have
$$
1 \leq N_L \leq n.
$$
\end{thm}

Notice that Polterovich \cite{Polterovich91'} gave an example of monotone Lagrangian submanifold
$L \subset \mathbb C^n$ with $N_L=n$, so these bounds are sharp.

Let us now turn to the case $(M,\omega)=({\mathbb C}P^n,\omega_{\rm FS})$.
We start with a preliminary lemma.

\begin{lem} \label{lem:minimal Maslov}
Let $L$ be a Lagrangian submanifold in ${\mathbb C}P^n$.
If $3 H_1(L;\mathbb Z)=0$, then $L$ is monotone and its minimal Maslov number is given by $N_L=2(n+1)/3$.
\end{lem}

\noindent
{\bf Proof.} \
Denote by $H_2^D$ the image of the Hurewicz homomorphism \\
$\pi_2({\mathbb C}P^n,L) \to H_2({\mathbb C}P^n,L;\mathbb Z)$.
Consider the following exact sequence
$$
\cdots \rightarrow H_2(L;\mathbb Z)
\stackrel{i_*}{\rightarrow} H_2({\mathbb C}P^n;\mathbb Z)
\stackrel{j_*}{\rightarrow} H_2({\mathbb C}P^n,L;\mathbb Z)
\stackrel{\partial_*}{\rightarrow} H_1(L;\mathbb Z) \rightarrow \cdots.
$$
For any $a \in H_2^D$, since $H_1(L;\mathbb Z)$ is $3$-torsion,
there exists $S \in H_2({\mathbb C}P^n; \mathbb Z)$ such that
$j_*(S)=3a \in H_2({\mathbb C}P^n,L; \mathbb Z)$.
Since $\mu(j_*(S))=2c_1(S)$, we have $\mu(a)=2 c_1(S)/3$, where
$c_1(S):=\langle c_1({\mathbb C}P^n),S \rangle \in \mathbb Z$
is the first Chern number of $S \in H_2({\mathbb C}P^n;\mathbb Z)$.
This equation deduces that $L$ is monotone, because ${\mathbb C}P^n$ is monotone as a symplectic manifold.
Let $S_0$ be the generator of $H_2({\mathbb C}P^n; \mathbb Z) \cong \mathbb Z$,
then we have $S=k S_0$ for some $k \in \mathbb Z$, and hence 
$$
\mu(a)=\frac{2}{3}k(n+1).
$$
This formula implies that $N_L \geq 2(n+1)/3$.
On the other hand, note that $1 \leq N_L \leq n+1$ holds for any monotone Lagrangian submanifold in ${\mathbb C}P^n$
from (\ref{eq:Maslovbound}).
Hence, we have $N_L=2(n+1)/3$.
\hfill \qed

\smallskip

Let $L \subset {\mathbb C}P^n$ be a Lagrangian submanifold which satisfies that $3H_1(L;\mathbb Z)=0$.
Since $N_L$ must be positive integer, by Lemma \ref{lem:minimal Maslov}, we have $3\,|\,n+1$ and $N_L$ is even.
It implies that $L$ is orientable \cite{Arnold67} (see also \cite[Proposition 3.1]{Damian12}).
Now we assume that $n=2$.
In this case $N_L=2$ holds.
The orientability of $L$ and the assumption that $H_1(L;\mathbb Z)$ is $3$-torsion implies that
$L \subset {\mathbb C}P^2$ is a Lagrangian sphere.
But, such $L$ cannot exist.
Indeed, if $L$ is simply connected, then we obtain $n_L=\# G_L(x)=1$, and
by (\ref{eq:Ono}) we have $N_L=6$, which is a contradiction.
Thus we have proved the first part of Theorem \ref{thm:main1}.

\smallskip

\noindent
{\bf Proof of Theorem \ref{thm:main1}\,(1).} \
By Lemma \ref{lem:minimal Maslov}, the assumption that $H_1(L;\mathbb Z)$ is $3$-torsion and $n=5$
imply that $L$ is monotone and $N_L=4$.
Consider the subcritical polarization $({\mathbb C}P^{n+1},\omega_{\rm FS},J,{\mathbb C}P^n)$ and
the Lagrangian circle bundle $S^1 \to \Gamma_L \stackrel{\pi}{\to} L$ constructed by Biran
\cite[Proposition 4.1.A]{Biran06}.
Then $\Gamma_L \subset \mathbb C^6$ is a monotone Lagrangian submanifold with the same minimal Maslov
number as $L$ ,i.e., $N_{\Gamma_L}=4$.
Since $\Gamma_L$ is displaceable in $\mathbb C^6$, we have $HF(\Gamma_L,\Gamma_L)=0$.

On the other hand, since $\nu=[(\dim \Gamma_L +1)/N_{\Gamma_L}]=[7/4]=1$,
Oh-Biran's spectral sequence $\{ E_r^{p,q},d_r \}$ collapses at stage $E_2$,
i.e., $E_2^{p,q}= \cdots =E_\infty^{p,q}$.
Here,
$$
E_2^{p,q}=
\frac{\ker \left( [\partial_1]:H^{p+q-pN_{\Gamma_L}}(\Gamma_L; \mathbb Z_2) \to H^{p+1+q-(p+1)N_{\Gamma_L}}(\Gamma_L; \mathbb Z_2) \right)}
{{\rm im} \left( [\partial_1]:H^{p-1+q-(p-1)N_{\Gamma_L}}(\Gamma_L; \mathbb Z_2) \to H^{p+q-pN_{\Gamma_L}}(\Gamma_L; \mathbb Z_2) \right)}
\otimes \Lambda^{pN_{\Gamma_L}}.
$$
Putting $p=0$ and $N_{\Gamma_L}=4$, we have
$$
E_2^{0,q}=
\frac{\ker \left( [\partial_1]:H^q(\Gamma_L; \mathbb Z_2) \to H^{q-3}(\Gamma_L; \mathbb Z_2) \right)}
{{\rm im} \left( [\partial_1]:H^{q+3}(\Gamma_L; \mathbb Z_2) \to H^q(\Gamma_L; \mathbb Z_2) \right)}.
$$
By Theorem \ref{thm:Oh-Biran}\,(5), we obtain the following exact sequence for any $q \in \mathbb Z$:
$$
H^{q+3}(\Gamma_L; \mathbb Z_2) \stackrel{[\partial_1]}{\longrightarrow} H^q(\Gamma_L; \mathbb Z_2)
\stackrel{[\partial_1]}{\longrightarrow} H^{q-3}(\Gamma_L; \mathbb Z_2).
$$
Note that $\Gamma_L \subset \mathbb C^6$ is orientable because $N_{\Gamma_L}$ is even (see \cite{Arnold67}).
Together with Poincar\'e duality and the universal coefficient theorem we have
\begin{eqnarray*}
& & H^3(\Gamma_L; \mathbb Z_2) \cong H^0(\Gamma_L; \mathbb Z_2) \oplus H^6(\Gamma_L; \mathbb Z_2)
 \cong \mathbb Z_2 \oplus \mathbb Z_2, \\
& & H^1(\Gamma_L; \mathbb Z_2) \cong H^4(\Gamma_L; \mathbb Z_2) \cong
 H^2(\Gamma_L; \mathbb Z_2) \cong H^5(\Gamma_L; \mathbb Z_2).
\end{eqnarray*}
The assumption for $H_1(L;\mathbb Z)$ means that it is the direct sum of some copies of ${\mathbb Z}_3$.
By the universal coefficient theorem, we obtain
\begin{eqnarray*}
& & H_1(L; \mathbb Z_2)
 \cong H_1(L; \mathbb Z) \otimes \mathbb Z_2 \oplus {\rm Tor}(H_0(L; \mathbb Z), \mathbb Z_2) \cong 0, \\
& & H^1(L; \mathbb Z_2)
 \cong {\rm Hom}(H_1(L; \mathbb Z), \mathbb Z_2) \oplus {\rm Ext}(H_0(L; \mathbb Z), \mathbb Z_2) \cong 0.
\end{eqnarray*}
By $\mathbb Z_2$-Poincar\'e duality, we have $H^4(L; \mathbb Z_2) \cong H_1(L; \mathbb Z_2) \cong 0$.

Now we use $\mathbb Z_2$-Gysin exact sequence for the fibration $S^1 \to \Gamma_L \stackrel{\pi}{\to} L$:
$$
0 \to H^1(\Gamma_L; \mathbb Z_2) \to H^0(L; \mathbb Z_2) \cong \mathbb Z_2 \stackrel{\cup \chi}{\to}
 H^2(L; \mathbb Z_2) \to H^2(\Gamma_L; \mathbb Z_2) \to 0,
$$
where $\cup \chi$ denotes the cup product by the (classical) $\mathbb Z_2$-Euler class $\chi$ of the fibration.
From it we can deduces that $H^1(\Gamma_L; \mathbb Z_2)$ is isomorphic to either $0$ or $\mathbb Z_2$.
But the case $H^1(\Gamma_L; \mathbb Z_2)=0$ cannot occur.
Indeed, we have $H^1(\Gamma_L; \mathbb R) \neq 0$ by Theorem \ref{thm:Gromov}
and it means that $H^1(\Gamma_L; \mathbb Z)$ contains at least one $\mathbb Z$-factor.
Then the universal coefficient theorem yields $H^1(\Gamma_L; \mathbb Z_2) \neq 0$.
Therefore, $H^1(\Gamma_L; \mathbb Z_2)$ must be isomorphic to $\mathbb Z_2$.
Then $H^4(\Gamma_L; \mathbb Z_2) \cong H^2(\Gamma_L; \mathbb Z_2) \cong H^5(\Gamma_L; \mathbb Z_2) \cong \mathbb Z_2$,
and hence the above exact sequence implies that
$H^2(L; \mathbb Z_2) \cong \mathbb Z_2$.
Again from $\mathbb Z_2$-Gysin exact sequence
$$
0 \to H^3(L; \mathbb Z_2) \to \mathbb Z_2 \oplus \mathbb Z_2 \to H^2(L; \mathbb Z_2) \cong \mathbb Z_2 \to 0,
$$
we have $H^3(L; \mathbb Z_2) \cong \mathbb Z_2$,
which concludes the proof.
\hfill \qed

\smallskip

In the above proof, we extracted informations about the cohomology of $L$ from the vanishing of
$HF(\Gamma_L,\Gamma_L)$.
This strategy was introduced by Biran in \cite{Biran06}.
Now we present a new approach.
If a Lagrangian submanifold $L$ itself satisfies that $HF(L,L)=0$,
then we can deduce from it directly somewhat stronger restrictions on the cohomology of $L$.

\begin{pro} \label{pro:rigidity}
Let $L$ be a Lagrangian submanifold in ${\mathbb C}P^n$ with $3 H_1(L;\mathbb Z)=0$.
Then $3\,|\,n+1$ and $N_L=2(n+1)/3$.
Moreover, if $HF(L,L)=0$ holds, then the cohomology of $L$ with coefficients in $\mathbb Z_2$ satisfies the following$:$
\begin{enumerate}
\item $H^0(L; \mathbb Z_2) \cong H^\frac{n+1}{3}(L; \mathbb Z_2) \cong H^\frac{2n-1}{3}(L; \mathbb Z_2) \cong H^n(L; \mathbb Z_2) \cong \mathbb Z_2$, \smallskip
\item $H^q(L; \mathbb Z_2) \cong 0$ \quad for \quad $q=\frac{n+4}{3},\frac{n+7}{3},\ldots,\frac{2n-4}{3}$, \smallskip
\item $H^q(L; \mathbb Z_2) \cong H^{q+\frac{2n-1}{3}}(L; \mathbb Z_2)$ \quad for \quad $q=1,2,\ldots,\frac{n-2}{3}$,
\end{enumerate}
where item $(2)$ holds when $n \geq 8$,
and the Euler characteristic $\chi(L)$ of $L$ is equal to zero.
\end{pro}

\noindent
{\bf Proof.} \
The first part has already proved in Lemma \ref{lem:minimal Maslov}.
For the second part, we use Oh-Biran's spectral sequence.
Assume that $HF(L,L)=0$.
Since $\nu=[(\dim L +1)/N_L]=[3/2]=1$,
the spectral sequence $\{ E_r^{0,q}, d_r \}$ collapses at stage $r=2$.
Recall that
$$
E_2^{0,q}=
\frac{\ker \left( [\partial_1]:H^q(L; \mathbb Z_2) \to H^{q+1-N_L}(L; \mathbb Z_2) \right)}
{{\rm im} \left( [\partial_1]:H^{q-1+N_L}(L; \mathbb Z_2) \to H^q(L; \mathbb Z_2) \right)}.
$$
Since $\oplus_{q \in \mathbb Z} E_2^{0,q} \cong HF(L,L)=0$ holds by the assumption and Theorem \ref{thm:Oh-Biran}\,(5),
we obtain the following exact sequence
\begin{equation} \label{eq:rigidity}
H^{q+\frac{2n-1}{3}}(L; \mathbb Z_2) \stackrel{[\partial_1]}{\longrightarrow} H^q(L; \mathbb Z_2)
\stackrel{[\partial_1]}{\longrightarrow} H^{q-\frac{2n-1}{3}}(L; \mathbb Z_2)
\end{equation}
for any $q \in \mathbb Z$.
From it we have the following isomorphisms
\begin{eqnarray*}
& &[\partial_1] : H^{\frac{2n-1}{3}}(L; \mathbb Z_2) \to H^0(L; \mathbb Z_2) \cong \mathbb Z_2, \\
& &[\partial_1] : H^n(L; \mathbb Z_2) \cong \mathbb Z_2 \to H^{\frac{n+1}{3}}(L; \mathbb Z_2),
\end{eqnarray*}
which yield item (1).
Similarly, items (2) and (3) are easily obtained from (\ref{eq:rigidity}).
For the last part, by item (2), we obtain
$$
\chi(L)=\sum_{q=0}^\frac{n+1}{3} (-1)^q \dim_{\mathbb Z_2}H^q(L; \mathbb Z_2)
 + \sum_{q=0}^\frac{n+1}{3} (-1)^{q+\frac{2n-1}{3}} \dim_{\mathbb Z_2}H^{q+\frac{2n-1}{3}}(L; \mathbb Z_2).
$$
Since $(2n-1)/3$ is odd, $\chi(L)=0$ holds from items (1) and (3).
\hfill \qed

\smallskip

To find examples of Lagrangian submanifolds of ${\mathbb C}P^n$ satisfying that $HF(L,L)=0$,
the following result and the next lemma are useful.

\begin{thm}[Biran-Khanevsky \cite{Biran-Khanevsky}, Corollary 1.1.2] \label{thm:B-K}
Suppose that $\Sigma$ appears as a symplectic hyperplane section in a symplectic manifold $M$
such that $W:=M \setminus \Sigma$ is subcritical.
Let $L \subset \Sigma$ be a monotone Lagrangian submanifold with $N_L \geq 3$.
Denote by $\mathcal{N} \to \Sigma$ the normal bundle of $\Sigma$ in $M$.
If $HF(L,L) \neq 0$, then the classical $\mathbb Z_2$-Euler class $e \in H^2(L; \mathbb Z_2)$
of the restriction $\mathcal{N}|_L$ is non-trivial.
\end{thm}

In the case that $\Sigma={\mathbb C}P^n \subset M={\mathbb C}P^{n+1}$ we can apply Theorem \ref{thm:B-K},
because $W \cong \mathbb C^{n+1}$ is subcritical.
Next, we give a criterion for $e(\mathcal{N}|_L)$ to be vanished.

\begin{lem} \label{lem:euler}
Let $L \subset {\mathbb C}P^n$ be a Lagrangian submanifold.
If $n$ is even and $L$ is orientable, then $e(\mathcal{N}|_L)=0$.
\end{lem}

\noindent
{\bf Proof.} \
Let us consider $\Sigma={\mathbb C}P^n$ as a linear hyperplane of $M={\mathbb C}P^{n+1}$ as above.
We take the standard complex structure $J$ on ${\mathbb C}P^{n+1}$. 
Then $\omega_{\rm FS}(\cdot,J\cdot)$ defines the Fubini-Study metric on ${\mathbb C}P^{n+1}$.
On ${\mathbb C}P^n$, we use the induced metric from it.
Then we have a decomposition
$$
T({\mathbb C}P^n)|_L \oplus \mathcal{N}|_L \cong T({\mathbb C}P^{n+1})|_L
$$
as complex vector bundles on $L$.
Denote by $c_1(E)$ the first Chern class of a complex vector bundle $E$.
Using the fact that $c_1(T({\mathbb C}P^{n+1}))=(n+2)\alpha$,
where $\alpha \in H^2({\mathbb C}P^{n+1}; \mathbb Z)$ is a generator,
we have
$$
c_1(T({\mathbb C}P^n)|_L) + c_1(\mathcal{N}|_L) = (n+2)\alpha|_L,
$$
it yields
$$
w_2(TL \oplus NL) + w_2(\mathcal{N}|_L) = 0 \quad (\mathrm{mod} \ 2),
$$
where $w_2$ is the second Stiefel-Whitney class.
Since $L \subset {\mathbb C}P^n$ is Lagrangian, the tangent bundle $TL$ of $L$
is isomorphic to the normal bundle $NL$ of $L$ in ${\mathbb C}P^n$.
Hence, we have
$$
w_2(TL \oplus NL) = 2w_2(TL) + w_1(TL) \cup w_1(TL) = 0 \quad (\mathrm{mod} \ 2),
$$
because $L$ is orientable.
Since the class $w_2$ of a two-dimensional real vector bundle equals the
$\mathbb Z_2$-Euler class $e$, we obtain
$e(\mathcal{N}|_L) = 0$. 
\hfill \qed

\smallskip

\noindent
{\bf Proof of Theorem \ref{thm:main1}\,(2).} \
By Lemma \ref{lem:minimal Maslov}, $L$ is monotone and $N_L=6$, and so $L$ is orientable.
Hence, the above lemma yields $e(\mathcal{N}|_L)=0$.
From Theorem \ref{thm:B-K} we have $HF(L,L)=0$,
and by Proposition \ref{pro:rigidity} we obtain
\begin{eqnarray*}
& & H^0(L; \mathbb Z_2) \cong H^3(L; \mathbb Z_2) \cong H^5(L; \mathbb Z_2) \cong H^8(L; \mathbb Z_2) \cong \mathbb Z_2, \\
& & H^4(L; \mathbb Z_2) \cong 0, \
 H^1(L; \mathbb Z_2) \cong H^6(L; \mathbb Z_2), \ H^2(L; \mathbb Z_2) \cong H^7(L; \mathbb Z_2).
\end{eqnarray*}
Since $H_1(L;\mathbb Z)$ is $3$-torsion, the universal coefficient theorem yields
$H_1(L;{\mathbb Z}_2) \cong 0$ and $H^1(L;{\mathbb Z}_2) \cong 0$.
Hence we have $H^6(L; \mathbb Z_2) \cong H^1(L; \mathbb Z_2) \cong 0$, and
$H^2(L; \mathbb Z_2) \cong H^7(L; \mathbb Z_2) \cong H_1(L; \mathbb Z_2) \cong 0$.
Therefore,
$$
H^q(L; \mathbb Z_2)=
\left\{
 \begin{array}{cl}
  \mathbb Z_2 & \quad (q=0,3,5,8)\\
  0 & \quad (\mathrm{otherwise}).
 \end{array}
\right.
$$
It is actually isomorphic to $H^*(SU(3)/\mathbb Z_3; \mathbb Z_2)$ as graded vector spaces.

Next we consider the product structure on $H^*(L;\mathbb Z_2)$.
Denote by $1, x_3, x_5$ and $x_8$ the generator of $H^q(L;\mathbb Z_2)$ for $q=0,3,5$ and $8$, respectively.
We first claim that $x_3 \star x_5=x_8$.
Observe that
$$
[\partial_1](x_8)=x_3,\ [\partial_1](x_5)=1,\ [\partial_1](x_3)=0 \in H^*(L;\mathbb Z_2).
$$
By Leibniz rule, we have
$$
[\partial_1](x_3 \star x_5) = [\partial_1](x_3) \cdot x_5 + x_3 \cdot [\partial_1](x_5)
 =x_3 \cdot 1 = x_3 = [\partial_1](x_8).
$$
Since $[\partial_1]: H^8(L;\mathbb Z_2) \to H^3(L;\mathbb Z_2)$ is an isomorphism (see Proposition \ref{pro:rigidity}),
we obtain $x_3 \star x_5=x_8$.
Similarly, we have
$$
[\partial_1](x_3 \star x_3) = 0,\ [\partial_1](x_5 \star x_5) = 2x_5 =0,
$$
which yield that $x_3 \star x_3=0$ and $x_5 \star x_5=0$.
Thus we complete the proof of (2) in Theorem \ref{thm:main1},
because on $V_1=H^*(L;\mathbb Z_2)$ the induced product from $\star$ coincides with the cup product on $H^*(L;\mathbb Z_2)$.
\hfill \qed

\smallskip

\noindent
{\bf Proof of Theorem \ref{thm:main1} (3).} \
The proof is the same as that of (2) above.
By Lemma \ref{lem:minimal Maslov}, $L$ is monotone and $N_L=18$, and hence $L$ is orientable.
Lemma \ref{lem:euler} yields $e(\mathcal{N}|_L)=0$, and hence we obtain $HF(L,L)=0$ from Theorem \ref{thm:B-K}.
Using the assumption that $H^i(L; \mathbb Z_2)=0 \ (i=2,3,4)$, as in the proof of (2) above, we obtain
$$
H^0(L; \mathbb Z_2) \cong H^9(L; \mathbb Z_2) \cong H^{17}(L; \mathbb Z_2) \cong H^{26}(L; \mathbb Z_2) \cong \mathbb Z_2,
$$
and for the other $q$ we have
$$
H^q(L; \mathbb Z_2) \cong 0
$$
except for $q=5,6,7,22,23,24$.

Since $\dim_{\mathbb Z_2} H^{26-q}(L; \mathbb Z_2)=\dim_{\mathbb Z_2} H^q(L; \mathbb Z_2)=0$ for $q=2,3,4$,
we obtain
$$
H^q(L; \mathbb Z_2) \cong H^{q+17}(L; \mathbb Z_2) \cong 0
$$
for $q=5,6,7$.
Therefore,
$$
H^q(L; \mathbb Z_2)=
\left\{
 \begin{array}{cl}
  \mathbb Z_2 & \quad (q=0,9,17,26)\\
  0 & \quad (\mathrm{otherwise}).
 \end{array}
\right.
$$
The proof of the isomorphism as algebras is completely same as that of (2) above,
so we omit it.
\hfill \qed

\smallskip

Finally, we shall prove that $\chi(L)=0$ for a Lagrangian submanifold $L \subset {\mathbb C}P^n$ with $3 H_1(L;\mathbb Z)=0$.
It suffices to consider the case where $n$ is even.
Since $L$ is orientable, monotone and $N_L \geq 6$, by Lemma \ref{lem:euler}, we have $e(\mathcal{N}|_L)=0$,
and hence $HF(L,L)=0$ by Theorem \ref{thm:B-K}.
Therefore, Proposition \ref{pro:rigidity} yields $\chi(L)=0$.
Thus we complete the proof of Theorem \ref{thm:main1}.

\section{Floer cohomology of a model Lagrangian submanifold}

In this section, we shall calculate the Floer cohomology of model Lagrangian submanifolds with coefficients in $\mathbb Z_2$,
especially the example $SU(p)/{\mathbb Z}_p \subset {\mathbb C}P^{p^2-1}$.

Combining Theorem \ref{thm:B-K} with Lemma \ref{lem:euler} in the previous section, we obtain the following vanishing result
of the Floer cohomology for orientable monotone Lagrangian submanifolds in ${\mathbb C}P^{2n}$.

\begin{pro} \label{pro:vanishHF}
Let $L \subset {\mathbb C}P^n$ be an orientable monotone Lagrangian submanifold with $N_L \geq 3$.
If $n$ is even, then the Floer cohomology $HF(L,L)$ with coefficients in $\mathbb Z_2$ vanishes.
\end{pro}

Notice that the assumption that $L$ is orientable and the dimension $n$ is even is necessary.
The real projective space ${\mathbb R}P^n \subset {\mathbb C}P^n$ well illustrates it.
Indeed, in this case $HF({\mathbb R}P^n,{\mathbb R}P^n) \neq 0$ and ${\mathbb R}P^n$ is orientable if and only if $n$ is odd.
We apply this proposition to the list in Proposition \ref{pro:A-O}.

\begin{cor} \label{cor:vanishingHF}
Let $L$ be one of the following Lagrangian submanifolds$:$
\begin{enumerate}
\item $L=SU(p)/{\mathbb Z}_p \subset {\mathbb C}P^{p^2-1}$, where $p \,(\geq 3)$ is an odd number,

\item $L=SU(p)/(SO(p) \mathbb Z_p) \subset {\mathbb C}P^{p(p+1)/2-1}$, where $p$ is prime and $4 \nmid p+1$,

\item $L=E_6/(F_4 \mathbb Z_3) \subset {\mathbb C}P^{26}$.
\end{enumerate}


\noindent
Then the Floer cohomology $HF(L,L)$ with $\mathbb Z_2$-coefficients vanishes.
\end{cor}

\noindent
{\bf Proof.}\ 
(1) By the assumption, $n=p^2-1$ is even.
The Lagrangian submanifold $L=SU(p)/{\mathbb Z}_p \subset {\mathbb C}P^n$ is monotone and orientable.
Moreover, we have $N_L \geq 2p$ from (\ref{eq:MMineq}).
Hence, by Proposition \ref{pro:vanishHF}, we obtain $HF(L,L)=0$.

(2) Since $p$ is prime, we have $N_L=p+1$, and hence the monotone Lagrangian submanifold 
$L=SU(p)/(SO(p) \mathbb Z_p) \subset {\mathbb C}P^n$ is orientable.
The additional assumption that $4 \nmid p+1$ implies that
$n+1=p(p+1)/2$ is odd, that is, $n$ is even.
Thus, we have $HF(L,L)=0$ from Proposition \ref{pro:vanishHF}.

(3) Since $L$ is monotone with $N_L=18$ and orientable, and the dimension of $L$ is even, we obtain $HF(L,L)=0$.
\hfill\qed

\smallskip

Here let us review the following useful terminologies introduced by Biran and Cornea \cite[Definition 1.2.1]{Biran-Cornea09}.

\begin{dfn} \rm 
A Lagrangian submanifold $L \subset (M,\omega)$ is said to be {\it narrow} if $HF(L,L)=0$,
and {\it wide} if there exists an isomorphism $HF(L,L) \cong H^*(L;\mathbb Z_2) \otimes \Lambda$
(in general not canonical), where $\Lambda={\mathbb Z}_2[T,T^{-1}]$.
\end{dfn}

Note that all known monotone Lagrangian submanifolds are either narrow or wide.
Corollary \ref{cor:vanishingHF} says that $SU(p)/{\mathbb Z}_p \subset {\mathbb C}P^{p^2-1}$
is narrow if $p$ is an odd number greater than or equal to three.
Now we are going to study the case where $p$ is even.
To pursue it, the following theorem is essential.

\begin{thm}[Biran-Cornea \cite{Biran-Cornea09}, Theorem 1.2.2] \label{thm:Biran-Cornea1}
Let $L^n \subset (M^{2n},\omega)$ be a monotone Lagrangian submanifold.
Assume that its singular cohomology $H^*(L; \mathbb Z_2)$ is generated as a ring with the cup product by
$H^{\leq l}(L; \mathbb Z_2)$.
\begin{enumerate}
\item If $N_L > l$, then $L$ is either wide or narrow.
Moreover, if $N_L > l+1$, then $L$ is wide.
\item In case $L$ is narrow, then $L$ is uniruled of order $k$ with $k=\max\{ l+1,n+1-N_L \}$ if $N_L < l+1$,
and $k=l+1$ if $N_L=l+1$.
\end{enumerate}
\end{thm}

We shall apply it to the monotone Lagrangian submanifold $L:=SU(p)/{\mathbb Z}_p \subset {\mathbb C}P^{p^2-1}$,
where $p$ is even.
The cohomology ring of $L$ with $\mathbb Z_2$-coefficients was calculated by Baum and Browder \cite{BaumBrowder65}.
For $p=2^r n' \ (2 \nmid n')$, we have
\begin{eqnarray} \label{eq:ring}
& & H^*\left(L;{\mathbb Z}_2\right)
 \cong \wedge_2 (x_1,x_3,\ldots,\widehat{x_{2\cdot2^r-1}},\ldots,x_{2p-1}) \otimes {\mathbb Z}_2[y]/(y^{2^r}),
\end{eqnarray}
where $\widehat{\ }$ denotes omission, $\deg y=2$, and with the additional relation that $y=x_1^2$ if $r=1$
(see \cite[Corollary 4.2]{BaumBrowder65}).

If $p=2^r$ for a natural number $r$, then by (\ref{eq:ring}) the singular cohomology of $L$ is generated as a ring by
$H^{\leq l}(L; \mathbb Z_2)$, where $l=2p-3$.
Hence, we have
$$
N_L \geq 2p > l+1,
$$
which implies that $L$ is wide from Theorem \ref{thm:Biran-Cornea1} (1).

If $p$ is even and not a power of $2$, then (\ref{eq:ring}) implies that
$H^*(L; \mathbb Z_2)$ is generated as a ring by $H^{\leq l}(L; \mathbb Z_2)$, where $l=2p-1$.
Therefore,
$$
N_L \geq 2p > l=2p-1
$$
holds, and hence $L$ is either wide or narrow from Theorem \ref{thm:Biran-Cornea1} (1).

Including the case where $p$ is odd, we have proved the following

\begin{thm} 
The Floer cohomology of the monotone Lagrangian submanifold 
$L=SU(p)/{\mathbb Z}_p \subset {\mathbb C}P^{p^2-1}$
with coefficients in $\mathbb Z_2$ is given as follows$:$
\begin{enumerate}
\item If $p$ is odd, then $HF(L,L)=0$.
\item If $p$ is a power of $2$, then $HF(L,L) \cong H^*(L;\mathbb Z_2) \otimes \Lambda$.
\item If $p$ is even and not a power of $2$,
then $HF(L,L)$ is isomorphic to either $H^*(L;\mathbb Z_2) \otimes \Lambda$ or $0$.
\end{enumerate}
\end{thm}

In particular, we obtain

\begin{cor}
Let $p$ be a power of $2$.
Then $L=SU(p)/{\mathbb Z}_p \subset {\mathbb C}P^{p^2-1}$ is non-displaceable.
Moreover, for any $\phi \in \mathrm{Ham}({\mathbb C}P^{p^2-1},\omega_{\rm FS})$ such that
$L$ and $\phi L$ intersect transversally, we have
$$
\#(L \cap \phi L) \geq p\cdot2^{p-1}.
$$
\end{cor}

On the other hand, in the case where $p$ is odd, the vanishing of the Floer cohomology of $L$ with
$\mathbb Z_2$-coefficients provides no information about its non-displaceability.
However, at least in the case where $p=3$, we can show the following

\begin{pro} \label{thm:main2}
The Floer cohomology of the monotone Lagrangian submanifold $SU(3)/{\mathbb Z}_3 \subset {\mathbb C}P^8$
with coefficients in $\mathbb Z$ is nonvanishing.
In particular, it is non-displaceable in ${\mathbb C}P^8$.
\end{pro}

\noindent
{\bf Proof.} \
Denote by $L$ the monotone Lagrangian submanifold $SU(3)/{\mathbb Z}_3$ of ${\mathbb C}P^8$.
Then $N_L=6$ and $L$ is orientable.
Notice that the cohomology ring of $SU(p^r)/{\mathbb Z}_{p^r} \ (r \in \mathbb N)$ with ${\mathbb Z}_p$-coefficients,
where $p$ is an odd prime, was calculated by A. Borel.
By \cite[p.\,309]{Borel54}, we have
$$
H^*\left(\frac{SU(3)}{\mathbb Z_3};{\mathbb Z}_3\right) \cong \wedge (x_1,x_3) \otimes {\mathbb Z}_3[y]/(y^3),
$$
where $\deg x_1=1, \deg x_3=3$ and $\deg y=2$.
This implies that $x_1 \wedge x_3$ and $y^2$ generates $H^4(L;{\mathbb Z}_3)$, and hence
\begin{eqnarray} \label{eq:H^4}
H^4(L;{\mathbb Z}_3) \cong {\mathbb Z}_3 \oplus {\mathbb Z}_3.
\end{eqnarray}

On the other hand, since $L$ is spin, the Floer cohomology is well-defined with $\mathbb Z$-coefficients \cite{FOOO09}
and Oh-Biran's spectral sequence works with $\mathbb Z$-coefficients (see \cite{Damian12}).

\smallskip

\noindent
{\bf Assume:} $HF(L,L;\mathbb Z)=0$.

\smallskip

Since $\nu=[(\dim L +1)/N_L]=1$ for $L$,
the spectral sequence $\{ E_r^{p,q}, d_r \}$ collapses at stage $r=2$, hence $E_2^{p,q}= \cdots =E_\infty^{p,q}$.
Recall that
$$
E_2^{0,q}=
\frac{\ker \left( [\partial_1]:H^q(L; \mathbb Z) \to H^{q+1-N_L}(L; \mathbb Z) \right)}
{{\rm im} \left( [\partial_1]:H^{q-1+N_L}(L; \mathbb Z) \to H^q(L; \mathbb Z) \right)}.
$$
Since for every $p \in \mathbb Z$, $\oplus_{q \in \mathbb Z} E_2^{p,q} \cong HF(L,L;\mathbb Z)=0$ holds,
we obtain the following exact sequences for $q \in \mathbb Z$:
$$
H^{q+5}(L; \mathbb Z) \stackrel{[\partial_1]}{\longrightarrow} H^q(L; \mathbb Z)
\stackrel{[\partial_1]}{\longrightarrow} H^{q-5}(L; \mathbb Z),
$$
which yield that
$$
H^4(L; \mathbb Z) \cong 0, \quad
H^5(L; \mathbb Z) \cong H^0(L; \mathbb Z) \cong \mathbb Z.
$$
Here we consider the cohomology of $L$ with coefficients in ${\mathbb Z}_3$.
By the universal coefficient theorem, we have
$$
H^4(L; {\mathbb Z}_3)
\cong H^4(L; \mathbb Z) \otimes {\mathbb Z}_3 \oplus {\rm Tor}(H^5(L; \mathbb Z),{\mathbb Z}_3) \cong 0,
$$
which contradicts to (\ref{eq:H^4}).
Therefore, we have $HF(L,L;\mathbb Z) \neq 0$, and $L$ is non-displaceable in ${\mathbb C}P^8$.
\hfill\qed

\smallskip

Thus we complete proofs of all the results in Section 1.2.
From the results in this section, the following natural question arises:

\smallskip
 
\noindent
{\bf Question.}
Is the monotone Lagrangian submanifold $SU(p)/{\mathbb Z}_p$ of ${\mathbb C}P^{p^2-1}$
non-displaceable for {\it any} $p \geq 2$?
More generally, are all Lagrangian submanifolds of ${\mathbb C}P^n$
in Proposition \ref{pro:A-O} non-displaceable?

\section{Uniruling of model Lagrangian submanifolds}

In this section, as a further application, we prove the existence of pseudo-holomorphic disc
with its boundary on a model Lagrangian submanifold $L$.
The result (Corollary \ref{cor:uniruling}) is an immediate consequence of Theorem \ref{thm:Biran-Cornea1} and
Corollary \ref{cor:vanishingHF} in the previous section.

\smallskip

\noindent
{\bf Proof of Corollary \ref{cor:uniruling}.} \
Let $L$ be the monotone Lagrangian submanifold $SU(p)/{\mathbb Z}_p \subset {\mathbb C}P^{p^2-1}$,
where $p$ is an odd prime.
Then we have $N_L=2p$.
By Corollary \ref{cor:vanishingHF}, we know that $SU(p)/{\mathbb Z}_p$ is narrow.
Furthermore, the singular cohomology $H^*(L; \mathbb Z_2)$ is generated as a ring by
$H^{\leq 2p-1}(L; \mathbb Z_2)$ from (\ref{eq:cohSU}).
Therefore, Theorem \ref{thm:Biran-Cornea1}\,(2) yields that
$L$ is uniruled of order $N_L=2p$.



\begin{rem} \rm 
The same argument as the above is also applicable to the following two cases:

$\bullet$ $L \subset {\mathbb C}P^8$ which satisfies the assumption of Theorem \ref{thm:main1},

$\bullet$ $L \subset {\mathbb C}P^{26}$ which satisfies the assumption of Theorem \ref{thm:main1} and \\
\hspace*{5mm} $H^i(L; \mathbb Z_2) \cong 0 \ (i=2,3,4)$. \\
Both $L$ are uniruling of order $N_L$.
\end{rem}

\section{Concluding remarks}

Finally, we notice that results about homological rigidity as in Theorem \ref{thm:Biran1} and Theorem \ref{thm:main1}
are useful in understanding the geography of monotone Lagrangian submanifolds of ${\mathbb C}P^n$ for small $n$.







As we pointed out in Remark \ref{rem:Biran1}, a monotone Lagrangian submanifold $L$ in ${\mathbb C}P^n$ with $N_L=n+1$
must be $\mathbb Z_2$-homological ${\mathbb R}P^n$.
It seems that a Lagrangian submanifold $L$ with relatively large $N_L$ tends to be
$\mathbb Z_2$-homologically rigid.
Consider a Lagrangian submanifold $L \subset {\mathbb C}P^8$.
The assumption $3 H_1(L;\mathbb Z)=0$ implies that $L \subset {\mathbb C}P^8$ is monotone and $N_L=6$.
Unfortunately, the argument in Theorem \ref{thm:main1}\,(2) does not imply that
a monotone Lagrangian submanifold $L \subset {\mathbb C}P^8$ with $N_L=6$ satisfies that
$H^*(L; \mathbb Z_2) \cong H^*\left(SU(3)/{\mathbb Z_3}; \mathbb Z_2\right)$,
because the proof essentially uses the assumption $3 H_1(L;\mathbb Z)=0$ to deduce that $H^1(L;\mathbb Z_2)=0$.
To overcome this difficulty Damian's lifted Floer homology theory \cite{Damian12} may be useful.

In contrast to that, for the case $N_L=2$,
at least ${\mathbb C}P^3$ admits Lagrangian submanifolds with distinct $\mathbb Z_2$-homological types, 
$T^3_{\rm clif}$ and $SO(3)/D_3$.
In \cite{Chiang04}, Chiang constructed a Lagrangian submanifold $SO(3)/D_3$ in ${\mathbb C}P^3$,
where $D_3$ is the dihedral group.
Its minimal Maslov number is calculated as follows.

\begin{lem}
The Lagrangian submanifold $L:=SO(3)/D_3 \subset {\mathbb C}P^3$ is monotone and $N_L=2$.
\end{lem}

\noindent
{\bf Proof.} \
Since the integral homology groups of $L$ are
$$
H_0(L; \mathbb Z) \cong \mathbb Z, \quad H_1(L; \mathbb Z) \cong \mathbb Z_4, \quad
H_2(L; \mathbb Z) \cong 0, \quad H_3(L; \mathbb Z) \cong \mathbb Z
$$
(see \cite[Section 3]{Chiang04}), $L$ is orientable, and so is spin because $\dim_{\mathbb R}L=3$.
Moreover, $L$ is monotone from a similar way to Lemma \ref{lem:minimal Maslov}.
Hence, $N_L$ is even, and we have $n_L N_L=8$ by (\ref{eq:Ono}).
Combining $N_L \leq 4$, we obtain $N_L=2$ or $4$.
If $N_L=4$, then we have $\pi_1(L) \cong \mathbb Z_2$ from Theorem \ref{thm:Damian1} and Remark \ref{rem:Biran1}.
It is impossible, so $N_L=2$ holds.
\hfill\qed

\smallskip

\begin{rem} \rm
Damian proved that a monotone oriented Lagrangian submanifold $L$ of ${\mathbb C}P^n$ which is aspherical
satisfies $N_L=2$ (see \cite[Theorem 1.6]{Damian12}).
However, since $SO(3)/D_3$ is spherical, the above lemma shows that the converse of the Damian's result does not hold.
\end{rem}

We also note that Biran and Cornea \cite[Section 6.4]{Biran-Cornea09} gave a method to construct Lagrangian submanifolds in
${\mathbb C}P^n$ with $N_L=2$.

\section*{Acknowledgements}

The author would like to thank Professor Paul Biran for valuable discussions on and his interest in the subject of this paper
which inspired me to obtain the principal idea of the present paper.
He also thank the FIM at ETH Zurich for providing an excellent working atmosphere and
great hospitality in March in 2013.
He would like to thank Professor Kaoru Ono, whose insightful comments on this work
improved the contents of Section 5.
The author was partly supported by the Grant-in-Aid for Young Scientists (B) (No.~24740049), JSPS.

\begin{flushleft}
H.~Iriyeh

{\sc School of Science and Technology for Future Life\\
Tokyo Denki University\\
5 Senju-Asahi-Cho, Adachi-Ku\\
Tokyo, 120-8551 Japan}

{\it e-mail} : {\tt hirie@im.dendai.ac.jp}
\end{flushleft}

\end{document}